\documentclass[12pt]{amsart}
\usepackage{amssymb,latexsym,amsthm}
\usepackage{amsmath,amsfonts,amssymb}
\usepackage{graphicx}
\usepackage{color}
\usepackage{mathrsfs}

\newcommand\A{\mathscr A}
\newcommand\BB{\mathscr B}
\newcommand\R{\mathbb R}

\newcommand\V{\mathbb V}
\newcommand\W{\mathbb W}

\newcommand{\0}{\mbox{\boldmath $0$}}
\newcommand{\x}{\mbox{\boldmath $x$}}
\newcommand{\y}{\mbox{\boldmath $y$}}
\newcommand{\z}{\mbox{\boldmath $z$}}
\newcommand{\PP}{\mbox{\boldmath $p$}}
\newcommand{\e}{\mbox{\boldmath $e$}}
\newcommand{\ab}{\mbox{\boldmath $a$}}
\newcommand{\bb}{\mbox{\boldmath $b$}}
\newcommand{\cb}{\mbox{\boldmath $c$}}
\newcommand{\uv}{\mbox{\boldmath $u$}}
\newcommand{\vv}{\mbox{\boldmath $v$}}

\newcommand\gyr{\operatorname {gyr}}
\newcommand\Aut{\operatorname {Aut}}

\newcommand{\un}{|\!\|}

\newtheorem{theorem}{Theorem}

\newtheorem{cor}[theorem]{Corollary}
\newtheorem{prop}[theorem]{Proposition}

\theoremstyle{definition}
\newtheorem{definition}[theorem]{Definition}
\newtheorem{example}[theorem]{Example}

\theoremstyle{remark}

\begin{document}
\author{Toshikazu Abe}
\address{Department of Information Engineering, Faculty of Engineering, Niigata University, Japan}
\email{abebin08@gmail.com}

\author{Osamu Hatori}
\address{Department of Mathematics, Faculty of Science, Niigata University, Japan}
\email{hatori@math.sc.niigata-u.ac.jp}

\title[A Note on the Paper]
{A Note on the Proof of Theorem 13 in the Paper "Generalized Gyrovector Spaces and a Mazur-Ulam theorem"}

\keywords{gyrovector spaces, a Mazur-Ulam theorem}

\subjclass[]{47B49,46L05,51M10}%


\begin{abstract}
We give a revision of the  proof of a Mazur-Ulam theorem for generalized gyrovector spaces given in \cite{abehatori}. 
\end{abstract}
\maketitle
\section{A proof of a Mazur-Ulam theorem for ggv's}
We introduce a notion of the generalized gyrovector space and give a Mazur-Ulam theorem for the generalized gyrovector space as Theorem 13 and Corollary 14 in \cite{abehatori}. The essential part of this Mazur-Ulam theorem is  Theorem 13, and it is exhibited without a required corection of the proof in the print version of \cite{abehatori} although the website version is revised. One may revise the proof easily, but for the sake of the convenience of the readers and the completeness of the proof we give a revision of the proof of Theorem 13. Notations and terminologies are due to \cite{abehatori}. The following is Theorem 13 in \cite{abehatori}.
\begin{theorem}\label{GMU}
Let $(G_1,\oplus_1,\otimes_1)$ and $(G_2,\oplus_2,\otimes_2)$ be GGV's with
$\varrho_1$ and $\varrho_2$ being gyrometrics of $G_1$ and $G_2$, respectively. Suppose that $T:G_1\to G_2$ is a gyrometric preserving surjection. 
Then $T$ preserves the gyromidpoints; 
\[\PP(T\ab,T\bb) =T\PP(\ab,\bb)\]
for any pair $\ab,\bb\in G_1$.
\end{theorem}
\begin{proof}
Let $\ab,\bb\in G_1$ and $\z$ be the gyromidpoint of $\ab$ and $\bb$. 
Let $W$ be the family of all bijective gyrometric preserving maps $S:G_1\rightarrow G_1$ keeping the points $\ab$ and $\bb$ fixed, 
and set 
\begin{equation}
\lambda = \sup \{f(\varrho(S\z,\z)) :S\in W\} \in [0,\infty],
\end{equation}
where $f$ is the bijection which satisfies {\rm{(F1)}} and {\rm{(F2)}} of Proposition 18 in \cite{abehatori}. Note that in \cite{abehatori} the function $f$ in Proposition 18 is denoted by $f:\|G\|\to {\mathbb R}$, but it reads as $f:\|\phi(G)\|\to {\mathbb R}$. 
For $S\in W$ we have $\varrho(S\z,\ab) =\varrho(S\z,S\ab) =\varrho(\z,\ab)$, hence
\begin{equation}
\varrho(S\z,\z)\leq\varrho(S\z,\ab)\oplus'_1\varrho(\ab,\z) = 2\otimes'_1\varrho(\ab,\z), 
\end{equation}
hence $f(\varrho(Sz,z))\le f(2\otimes'_1\varrho(a,z))=2f(\varrho(a,z))$ which yields $\lambda < \infty$.

Let $\psi (\x)=2\otimes_1 \z\ominus_1 \x$ on $G_1$. If $S \in W$, then so also is $S^{\ast}=\psi S^{-1}\psi S$, 
and therefore $f(\varrho(S^{\ast}\z, \z))\leq\lambda$. 
Since $S^{-1}$ is a gyrometric preserving map, (p5) of Proposition 16 and Proposition 18 in \cite{abehatori} together imply that 
\begin{equation}
\begin{split}
\lambda \geq f(\varrho(S^{\ast}\z,\z)) &= f(\varrho(\psi S^{-1}\psi S\z, \z)) \\
                                     &= f(\varrho(S^{-1}\psi S\z, \z)) \\
                                     &= f(\varrho(\psi S\z, S\z)) \\
                                     &= f(2\otimes'_1 \varrho(S\z,\z)) \\
                                     &= 2f(\varrho(S\z,\z))
\end{split}
\end{equation}
for all $S\in W$, showing that $\lambda \geq 2 \lambda$. Thus $\lambda =0$, which means that $S\z=\z$ for all $S\in W$ since $f$ is linear and bijective. 

Let $T:G_1\rightarrow G_2$ be a bijective gyrometric preserving map. 
Let $\z '$ be the gyromidpoint of $T(\ab )$ and $T(\bb )$. To prove the theorem we must show that $T(\z )=\z '$. 
Let $\psi '(\y )=2\otimes_2 \z '\ominus_2 \y$ on $G_2$. 
Then the map $\psi T^{-1} \psi ' T$ is in $W$, whence $\psi T^{-1} \psi ' T(\z )=\z$. 
This implies that $\psi ' (T(\z ))=T(\z )$. 
Due to (p3) of  Proposition 16, we obtain $T(\z )=\z '$. 
\end{proof}


\end{document}